\documentclass{cqjm}
\usepackage{amssymb} \usepackage{mathtools}
\usepackage{booktabs,amsmath}
\usepackage{graphicx} 
\usepackage{cite}
\usepackage{enumerate}
\usepackage[expansion=true,protrusion=false]{microtype}
\usepackage{dsfont}				
\renewcommand{\theequation}{\arabic{section}.\arabic{equation}}

\selectfont

\usepackage{placeins}

\usepackage{xcolor}

\DeclareMathOperator{\F}		{\mathcal F}	
\DeclareMathOperator{\G}		{\Gamma}		

\allowdisplaybreaks

\begin{document}

\title{A Novel Property of \\Generalized Fibonacci Sequence in Grids}        

\sfootnotetext[]{\noindent\textbf{Date}: 2024-3-7\\  
\textbf{Foundation item}: Supported by Shaanxi Fundamental Science Research Project for Mathematics and Physics(22JSQ001) \& Shaanxi
Province College Student Innovation and Entrepreneurship Training Program(S202210699481, S202310699324X) \\
\textbf{Biographies}:
Zixian Yang (2002-), male, native of Mianyang, Sichuan, an undergraduate of Northwestern Polytechnical University, Zixian\_Yang@outlook.com;       
Jianchao Bai (1987-), male, native of Shangluo, Shaanxi, an
 associate professor of Northwestern Polytechnical University, engages in Numerical Algebra and Optimization.
}

\enauthor{Zixian Yang, Jianchao Bai}
{(~{\it School of Mathematics and Statistics, Northwestern Polytechnical University, Xi'an 710129, China}~)}

\runhead{A Novel Property of Generalized Fibonacci Sequence in Grids}

\pagestyle{myheadings}

\begin{abstract}     
  Fibonacci sequence, generated by summing the preceding two terms, is a classical sequence renowned for its elegant properties. In this paper, leveraging properties of generalized Fibonacci sequences and formulas for consecutive sums of equidistant subsequences, we investigate the ratio of the sum of numbers along main-diagonal and sub-diagonal of odd-order grids containing generalized Fibonacci sequences. We show that this ratio is solely dependent on the order of the grid, providing a concise and splendid identity.
\end{abstract}

\MR{11B39}  
\CLC{O156.4}{A}

\maketitle
\section{Introduction}
\setcounter{equation}{0}
The classical Fibonacci sequence, denoted by $\{F_{n}\}$,  has certain relationship with Lucas  sequence \cite{cite1,cite2,cite3}. It appears everywhere in nature and art \cite{sinhaFibonacciNumbersIts2017,schreiberFibonacciSequenceArt2022}, such as leaf arrangement in plants, pattern of the florets of a flower, as so on. Generally speaking, it   is  defined by
\[
F_{n}=F_{n-1}+F_{n-2}, \qquad \textrm{for}~ n\geq 2,
 \]
with $ F_{0}=0$ and $ F_{1}=1$. In 1843, Binet gave an explicit formula  $F_{n}=\frac{1}{\sqrt{5}}\left(a^{n}-b^{n}\right)$, where $a=\frac{1+\sqrt{5}}{2},b=\frac{1-\sqrt{5}}{2}$. As a variation, Lucas sequence,  defined by $L_{n}=L_{n-1}+L_{n-2}$ with $ L_{0}=2$ and $ L_{1}=1$, also has an explicit formula $L_{n}=a^n+b^n$. Both of these sequences has been reviewed in  \cite[Chapter 7]{krizekFibonacciLucasNumbers2021}. In recent years, Fibonacci sequence is applied in optimization algorithms \cite{chenComplexitySequentiallyLifting2021,tran-ngocPromisingApproachUsing2023} and mechanical structure design \cite{mengRevisitingFibonacciSpiral2023}, due to  its remarkable structural elegance and mathematical beauty. Modern science has shown a strong interest in the theory and application of the Golden Section and Fibonacci numbers as well \cite{APGeneralized}.\par

It is fully acknowledged that Fibonacci sequence satisfies  some beautiful identities, such as Robert Simson identity \cite{RSMFibonacci}, Gelin-Cesaro identity \cite{LEDHistory}, Hoggatt and Bergum identities \cite{HBProblem}, Zvonko Cerin identities \cite{ZCFactors} and so on.
Some of these identities mentioned are listed below:
\begin{itemize}
  \item Robert Simson identity:
$
  F_{n-1}F_{n+1}-F_{n}^{2}=(-1)^{n}.
$
  \item 
  Gelin-Cesaro identity:
$
  F_n^4-F_{n-2}F_{n-1}F_{n+1}F_{n+2}=1.
$
  \item 
  Hoggatt and Bergum identities:
\[
  \begin{aligned}
    F_{n}F_{n+3}^{2}-F_{n+2}^{3}&=(-1)^{n+1}F_{n+1}, \\
    F_{n+3}F_{n}^{2}-F_{n+1}^{3}&=\left(-1\right)^{n+1}F_{n+2}, \\
    F_{n}F_{n+3}^{2}-F_{n+4}F_{n+1}^{2}&=(-1)^{n+1}L_{n+2}, \\
    F_{n}L_{n+3}^{2}-F_{n+4}L_{n+1}^{2}&=(-1)^{n+1}L_{n+2}.
  \end{aligned}
  \] 
   \item  
   Zvonko Cerin identities:
\[
    \begin{aligned}
      \sum_{j=0}^{4i+3}F_{k+j}=F_{2i+2}L_{k+2i+3},\quad& \quad\sum_{j=0}^{4i+3}(-1)^jF_{k+j}=F_{2i+2}L_{k+2i},\\
      \sum_{j=0}^{4i+1}F_{k+j}=L_{2i+1}F_{k+2i+2},\quad& \quad\sum_{j=0}^{4i+1}(-1)^jF_{k+j}=L_{2i+1}F_{k+2i-1},\\
      \sum_{j=0}^{4i}F_{k+j}=F_{2i}L_{k+2i}+L_{2i+1}F_{k+2i},\quad& \quad\sum_{j=0}^{4i}(-1)^jF_{k+j}=F_{k+2i}L_{2i+1}-L_{k+2i}F_{2i}.
    \end{aligned}
    \]
\end{itemize}
 
Significantly, Zvonko Cerin identities reveal the fact that the summation of a part of Fibonacci sequence contains a specific factor in Fibonacci form. In this paper, we will introduce and prove a novel identity of generalized Fibonacci sequence
\begin{equation}
  \frac{\sum_{k=0}^{2n}G_{2k(n+1)}}{\sum_{k=0}^{2n}G_{2n(k+1)}}=c(n)
  \nonumber
\end{equation}
where $c(n)$ is a ratio solely dependent on $n$. 
This equation also reveals that the summation has a specific factor but is presented in a structurally symmetrical form of a division.
Additionally, we introduce a geometric interpretation for the equation above to provide   a simple yet intuitive  understanding.
\section{Preliminaries}
\setcounter{equation}{0}
\begin{definition}
  The generalized Fibonacci sequence $\{G_n\}_{n=-\infty}^\infty$ is defined by
  \[
    G_{n}=G_{n-1}+G_{n-2},\quad \textrm{with}~ G_{0}=A, G_{1}=B,
  \]
  where $A,B\in\mathbb{Z}$ such that $A^{2}+B^{2}\neq0$.
\end{definition}

Actually, the  condition $A,B\in\mathbb{R}$ such that $A^{2}+B^{2}\neq0$ is also workable, but we  only focus on the case of integers in this article. It is evident that  both   Fibonacci sequence and Lucas sequence are special cases of the above generalized Fibonacci sequence.

\begin{lemma}\label{lm1}
The generalized Fibonacci sequence $\{G_n\}_{n=-\infty}^{\infty}$ can be linearly expressed by Fibonacci sequence $\{F_n\}_{n=-\infty}^{\infty}$ and Lucas sequence $\{L_n\}_{n=-\infty}^{\infty}$, that is,
	\begin{equation}\label{eqlm1}
		G_{n}=\frac{1}{2}\bigl(-A+2B\bigr)F_{n}+\frac{1}{2}A L_{n}.
	\end{equation}
\end{lemma}\\
{\bf Proof.}
Under the rules of recursive formulas, $\{G_n\}$ can be uniquely determined by the initial  values  $G_0=A, G_1=B$. So can $\{F_n\}$ and $\{L_n\}$. Due to the determinant
\[
	\begin{vmatrix}F_{0}&L_{0}\\F_{1}&L_{1}\end{vmatrix}=\begin{vmatrix}0&2\\1&1\end{vmatrix}=-2\neq0,
\]
$G_0$ and $G_1$ can be linearly expressed by $F_0, F_1, L_0, L_1$. Let $G_{n}=k_{F}F_{n}+k_{L}L_{n}$, then Equation \eqref{eqlm1} can be proved by solving the following linear equation
\[
\begin{pmatrix}F_{0}&L_{0}\\F_{1}&L_{1}\end{pmatrix}\begin{pmatrix}k_{F}\\k_{L}\end{pmatrix}=\begin{pmatrix}G_{0}\\G_{1}\end{pmatrix}
\]
{
with respect to $(k_{F}, k_{L})^{\top}$. Obviously, its solution  is $(k_{F}, k_{L})^{\top}=\big(\frac{-A+2B}{2}, \frac{A}{2}\big)^{\top}$.} $\hfill\blacksquare$

\begin{lemma}\label{lm2}
	The sum of a continuous $n$ terms of the generalized Fibonacci subsequence $\{G_n\}_{n=-\infty}^{\infty}$ with distance $m \in\mathbb{N}$ has a formula:
  \begin{equation}\label{eqlm2}
    \sum_{k=n_1}^{n_2}G_{mk+i}=\frac1{a^m-b^m}\left[\left(G_{m+i}-b^{m}G_{i}\right)\frac{a^{mn_1}-a^{m(n_2+1)}}{1-a^{m}}-\left(G_{m+i}-a^{m}G_{i}\right)\frac{b^{mn_1}-b^{m(n_2+1)}}{1-b^m}\right],
  \end{equation}
where $i=0,1,2,\cdots,m-1$, $a=\frac{1+\sqrt{5}}{2}$ and $b=\frac{1-\sqrt{5}}{2}$.
\end{lemma}\\
{\bf Proof.}
Reference \cite{wangImprovingSumFormula2009} provides a concise proof for $\{G_{n}\}_{n=0}^{\infty}$. Since the first two terms of the sequence can be any initial values which are not all zero, it is true that $\{G_{n}\}_{n=-\infty}^{\infty}$ satisfies Equation \eqref{eqlm2}. $\hfill\blacksquare$

\begin{lemma}\label{lm3}
  For any even number $n$, Fibonacci sequence and Lucas sequence satisfy
  \[
    F_{m}L_{n}=F_{m+n}+F_{m-n}.
  \]
\end{lemma}\\
{\bf Proof.}
Because of the known equality $F_{m}L_{n}=F_{m+n}+\left(-1\right)^{n}F_{m-n}$, Lemma \ref{lm3} holds when $n$ is an even number. $\hfill\blacksquare$
\section{Property of generalized Fibonacci sequence}
\setcounter{equation}{0}
In this section, we provide a new discovery, shown in Theorem \ref{tr2}, for the generalized Fibonacci sequences within grids of any odd order.

\begin{figure}[h]
\def\figurename{Fig.}
  \centering
  \includegraphics[width=0.95\textwidth]{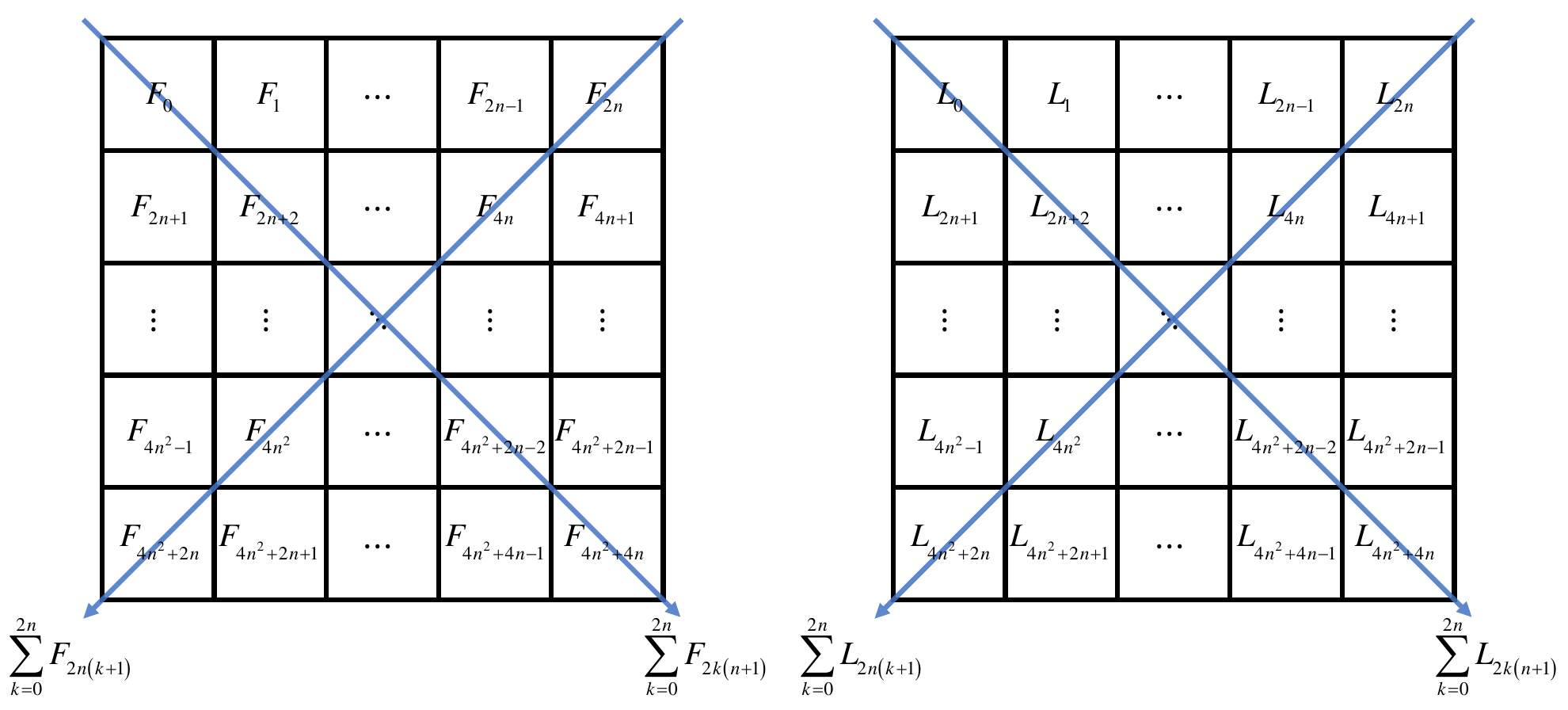}
  \caption{Fibonacci numbers and Lucas numbers arranged in odd order grids}
  \label{Fibonacci numbers and Lucas numbers arranged in odd order grids}
\end{figure}

\begin{theorem}\label{tr2}
  In any grid of odd order $2n+1$, the ratio between the sum of numbers along the main diagonal $\sum_{k=0}^{2 n} G_{2 k(n+1)}$ and the sum of numbers along the sub-diagonal $\sum_{k=0}^{2n}G_{2n(k+1)}$ of the generalized Fibonacci sequence $\{G_n\}$:
  \begin{equation}\label{eqtr2}
    \frac{\sum_{k=0}^{2n}G_{2k(n+1)}}{\sum_{k=0}^{2n}G_{2n(k+1)}}=c(n)
  \end{equation}
  is solely dependent on the value of $n$.
\end{theorem}\\
{\bf Proof.}
Based on  Lemma \ref{lm1}, if
\begin{equation}\label{eqFL}
  c_F(n):=\frac{\sum_{k=0}^{2 n} F_{2 k(n+1)}}{\sum_{k=0}^{2 n} F_{2 n(k+1)}}=\frac{\sum_{k=0}^{2 n} L_{2 k(n+1)}}{\sum_{k=0}^{2 n} L_{2 n(k+1)}}:=c_L(n)
\end{equation}
is true, Theorem \ref{tr2} can be proven.
As illustrated in Figures \ref{Fibonacci numbers and Lucas numbers arranged in odd order grids}, Fibonacci sequence and Lucas sequence are sequentially filled into square grids of order $2n+1$ starting from the 0th term. Thus our focus lies in demonstrating that the ratio between the sums of numbers along the main and sub-diagonals of Lucas sequence and Fibonacci sequence within a certain odd order grid remains consistent.

By Lemma \ref{lm2}, the two sides of the fractional line in $c_F(n)$ follows that
\begin{flalign}
  \sum_{k=0}^{2n}F_{2k(n+1)}=\sum_{k=0}^{2n}F_{2(n+1)k}&=\frac1{a^{2(n+1)}-b^{2(n+1)}}\left[\left(F_{2(n+1)}-b^{2(n+1)}F_0\right)\frac{a^{2(n+1)\cdot0}-a^{2(n+1)(2n+1)}}{1-a^{2(n+1)}}\right. \nonumber  \\ & \left.  \quad-\left(F_{2(n+1)}-a^{2(n+1)}F_0\right)\frac{b^{2(n+1)
  \cdot0}-b^{2(n+1)(2n+1)}}{1-b^{2(n+1)}}\right],\label{bjjc-1}\\
  \sum_{k=0}^{2n}F_{2n(k+1)}=\sum_{k=0}^{2n}F_{2nk+2n}&=\frac1{a^{2n}-b^{2n}}\left[\Big(F_{2n+2n}-b^{2n}F_{2n}\Big)\frac{a^{2n\cdot0}-a^{2n(2n+1)}}{1-a^{2n}} \right. \nonumber  \\ & \left. \quad -\Bigl(F_{2n+2n}-a^{2n}F_{2n}\Bigr)\frac{b^{2n\cdot0}-b^{2n(2n+1)}}{1-b^{2n}}\right],\label{bjjc-2}
\end{flalign}
where $a=\frac{1+\sqrt{5}}{2},b=\frac{1-\sqrt{5}}{2}$. By letting
\begin{equation}\label{10}
  \left\{
    \begin{aligned}
      d_F\left(n\right)&=-\frac{b^{2+2n}\left(-1+b^{4n(1+n)}\right)+a^{2+6n+4n^2}
      \left(-1+b^{2+2n}\right)-a^{2+2n}\left(-1+b^{2+6n+4n^2}\right)}{b^{2n}
      -b^{4n(1+n)}-a^{4n(1+n)}\left(-1+b^{2n}\right)+a^{2n}\left(-1+b^{4n(1+n)}\right)},\\
      d_L\left(n\right)&=\frac{2-b^{2+2n}+\left(-1+a^{2+2n}\right)b^{2+6n+4n^2}+a^{2+2n}
      \left[-1+a^{4n(1+n)}\left(-1+b^{2+2n}\right)\right]}{b^{2n}-b^{4n(1+n)}+a^{4n(1+n)}
      \left(-1+b^{2n}\right)+a^{2n}\left(1-2b^{2n}+b^{4n(1+n)}\right)},
    \end{aligned}
  \right.
\end{equation}
and combining the above equations  (\ref{bjjc-1}) and (\ref{bjjc-2}) together with  $F_{n}=\frac{1}{\sqrt{5}}\left(a^{n}-b^{n}\right)$, we can obtain
\[
  \begin{aligned}
      c_{F}(n)=\frac{\left(-1+a^{2n}\right)\left(-1+b^{2n}\right)}{\left(-1+a^{2+2n}\right)\left(-1+b^{2+2n}\right)}d_F(n).
  \end{aligned}
\]
Similarly, by Lemma \ref{lm2} and $L_{n}=a^{n}+b^{n}$, we have
\[
  \begin{aligned}
      c_{L}(n)=\frac{\left(-1+a^{2n}\right)\left(-1+b^{2n}\right)}{\left(-1+a^{2+2n}\right)\left(-1+b^{2+2n}\right)}d_L(n).
  \end{aligned}
\]

Based on the above discussions, to prove (\ref{eqFL}), we only need to prove $d_{F}(n)=d_{L}(n)$. Rearrange Equation \eqref{10} to have
\begin{equation}
  \begin{split}
  &\left\{
      \begin{aligned}
      d_F(n)&=\frac{-F_{2n+2}-F_{4n^2+4n}+F_{4n^2+6n+2}}{-F_{2n}-F_{4n^2+2n}+F_{4n^2+4n}} &=\frac{F_0-F_{2n+2}-F_{4n^2+4n}+F_{4n^2+6n+2}}{F_0-F_{2n}-F_{4n^2+2n}+F_{4n^2+4n}},\\
      d_L(n)&=\frac{2-L_{2n+2}+L_{4n^2+4n}-L_{4n^2+6n+2}}
      {-2+L_{2n}+L_{4n^2+2n}-L_{4n^2+4n}} &=\frac{L_0-L_{2n+2}+L_{4n^2+4n}-L_{4n^2+6n+2}}
      {-L_0+L_{2n}+L_{4n^2+2n}-L_{4n^2+4n}},
    \end{aligned}
  \right. \\
  \Rightarrow
  &\left\{
    \begin{aligned}
      d_F(n)&=\frac{-\left(a^{2n+2}-b^{2n+2}\right)-\left(a^{4n^2+4n}-b^{4n^2+4n}\right)+\left(a^{4n^2+6n+2}-b^{4n^2+6n+2}\right)}{-\left(a^{2n}-b^{2n}\right)-\left(a^{4n^2+2n}-b^{4n^2+2n}\right)+\left(a^{4n^2+4n}-b^{4n^2+4n}\right)},\\
      d_L(n)&=\frac{2-\left(a^{2n+2}+b^{2n+2}\right)+\left(a^{4n^2+4n}
      +b^{4n^2+4n}\right)-\left(a^{4n^2+6n+2}+b^{4n^2+6n+2}\right)}{-2+\left(a^{2n}
      +b^{2n}\right)+\left(a^{4n^2+2n}+b^{4n^2+2n}\right)-
      \left(a^{4n^2+4n}+b^{4n^2+4n}\right)}.
    \end{aligned}
  \right. \nonumber
  \end{split}
\end{equation}
So, it  is equivalent to proving
\[
  \begin{aligned}&\left(-F_{2n+2}-F_{4n^{2}+4n}+F_{4n^{2}+6n+2}\right)\left(-2+L_{2n}+
  L_{4n^{2}+2n}-L_{4n^{2}+4n}\right)\\=
  &\left(-F_{2n}-F_{4n^{2}+2n}+F_{4n^{2}+4n}\right)\left(2-L_{2n+2}+L_{4n^{2}+4n}-
  L_{4n^{2}+6n+2}\right),
  \end{aligned}
\]
that is,
\begin{equation}\label{bjc-3}
  \begin{aligned}
  &2F_{2+2n}+2F_{4n+4n^2}-2F_{2+6n+4n^2}-F_{2+2n}L_{2n}-F_{4n+4n^2}L_{2n} \\
  &+F_{2+6n+4n^2}L_{2n}-F_{2+2n}L_{2n+4n^2}-F_{4n+4n^2}L_{2n+4n^2}+F_{2+6n+4n^2}L_{2n+4n^2} \\
  &+F_{2+2n}L_{4n+4n^2}+F_{4n+4n^2}L_{4n+4n^2}-F_{2+6n+4n^2}L_{4n+4n^2} \\
  =&-2F_{2n}-2F_{2n+4n^2}+2F_{4n+4n^2}+F_{2n}L_{2+2n}+F_{2n+4n^2}L_{2+2n} \\
  &-F_{4n+4n^2}L_{2+2n}-F_{2n}L_{4n+4n^2}-F_{2n+4n^2}L_{4n+4n^2}+F_{4n+4n^2}L_{4n+4n^2} \\
  &+F_{2n}L_{2+6n+4n^2}+F_{2n+4n^2}L_{2+6n+4n^2}-F_{4n+4n^2}L_{2+6n+4n^2}.
  \end{aligned}
\end{equation}
Now, applying Lemma \ref{lm3}, we have
\[
  \begin{aligned}
  &\text{Left-hand-side  of (\ref{bjc-3})}\\
  =& 2F_{4n^2+4n}-F_2-F_{2n}+F_{2n+2}+F_{4n^2-2}-F_{4n^2+2n}-F_{4n^2+6n}-F_{8n^2+6n}  \\
  &+F_{8n^2+8n}-F_{4n^2+2n-2}-F_{4n^2+6n+2}+F_{4n^2+8n+2}+F_{8n^2+8n+2}-F_{8n^2+10n+2} \\
  =&\text{Right-hand-side of (\ref{bjc-3})}
  \end{aligned}
\]
Finally, Equation \eqref{eqtr2} is proven, and the   conclusion  is established. $\hfill\blacksquare$
\section{Examples}
\setcounter{equation}{0}
In this section, we  provide some examples of the property above to show the validity intuitively, and we shall present the proof in the case of order-3 grid.
\subsection{Fibonacci sequence in certain order girds}
\begin{figure}[htbp]
  \def\figurename{Fig.}
    \begin{minipage}{1\textwidth}
      \centering
      \resizebox{0.8\textwidth}{!}{\includegraphics{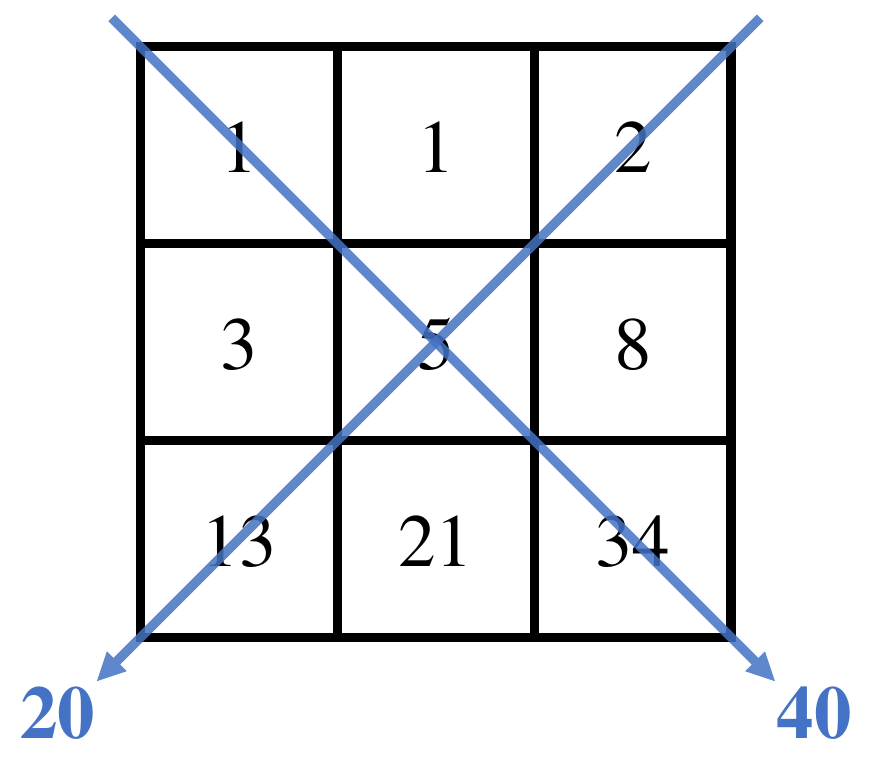}
  \hspace{5cm}\includegraphics{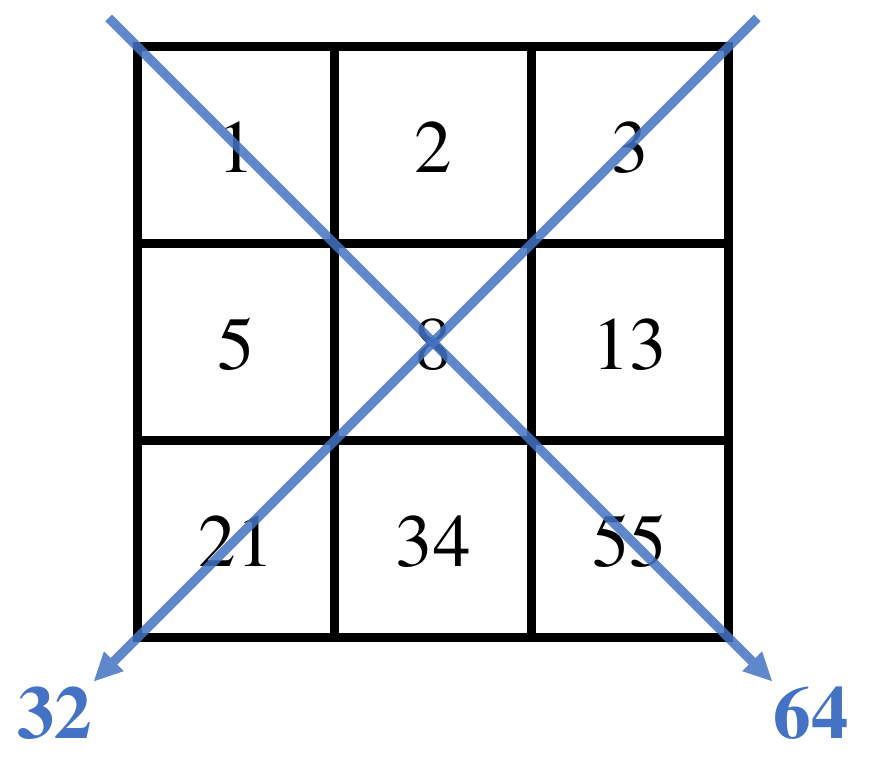}}
      \caption{Fibonacci numbers arranged in a $3\times 3$ grids} \label{bjcy-4}
    \end{minipage}
  \end{figure}

Given a $3\times 3$ grid, fill the 1st to 9th Fibonacci numbers into the grid in turn, and calculate the sum of the numbers on the main diagonal and the sub-diagonal, as shown in the left of Figure \ref{bjcy-4}. It is also observed that the 2nd to 10th Fibonacci numbers are operated similarly in the right of Figure \ref{bjcy-4}.

\begin{figure}[htbp]
  \def\figurename{Fig.}
    \centering
    \includegraphics[width=0.95\textwidth]{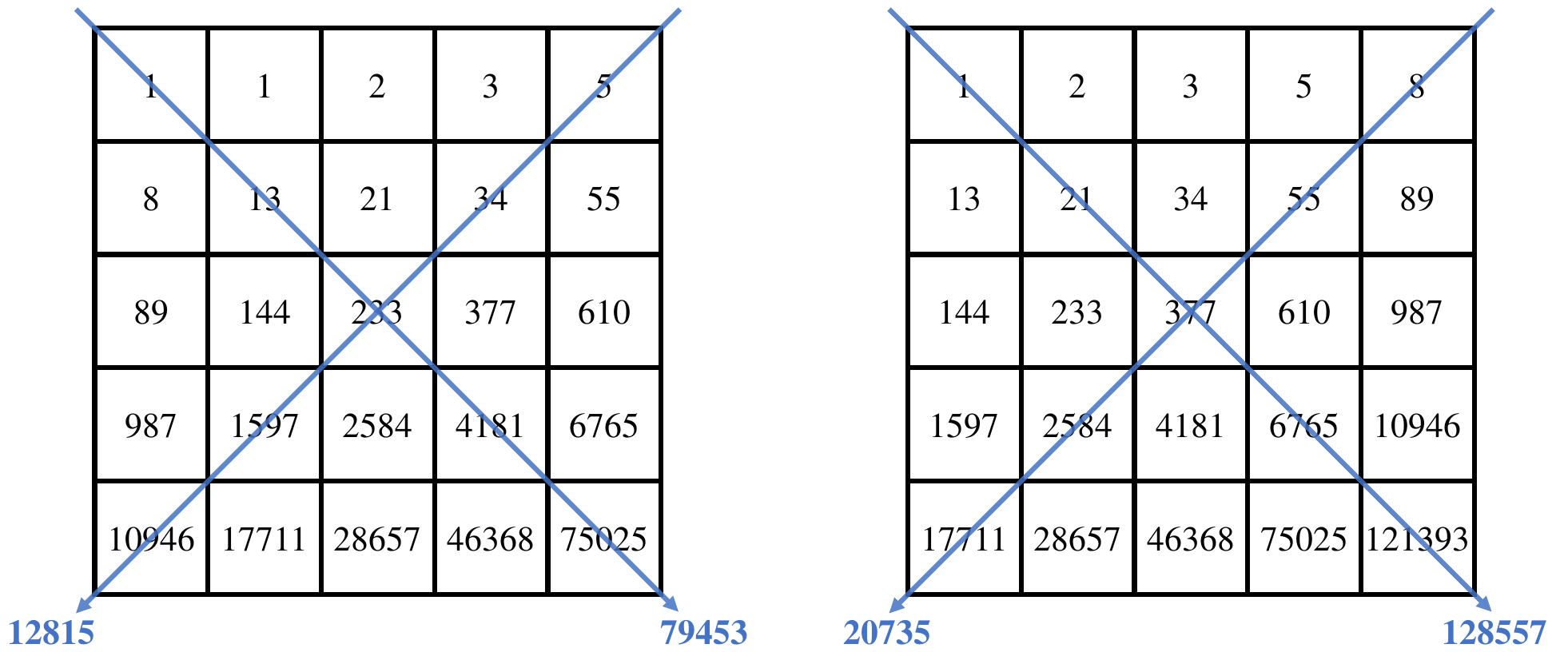}
    \caption{Fibonacci numbers arranged in a $5\times 5$ grid}
    \label{bjcy-5*5}
  \end{figure}

We can observe that $\frac{40}{20}=\frac{64}{32}=2$, that is, the ratio of the sum of the numbers on the main diagonal to the sum of the numbers on the sub-diagonal in the first of  Figure \ref{bjcy-4} is equal to that in the second of Figure \ref{bjcy-4}.
This observation is also true for the $5\times 5$ grid,
  that $\frac{79453}{12815} = \frac{128557}{20735}=6.2$ is a constant, as shown in Figure \ref{bjcy-5*5}. Similarly, it can be verified for $7\times 7$, $9\times 9$ and $11\times 11$, but not for $2\times 2$ and $4\times 4$.

\subsection{Generalized Fibonacci sequence in order-3 grid}
\begin{figure}[!htbp]
  \def\figurename{Fig.}
    \centering
    \includegraphics[width=0.45\textwidth]{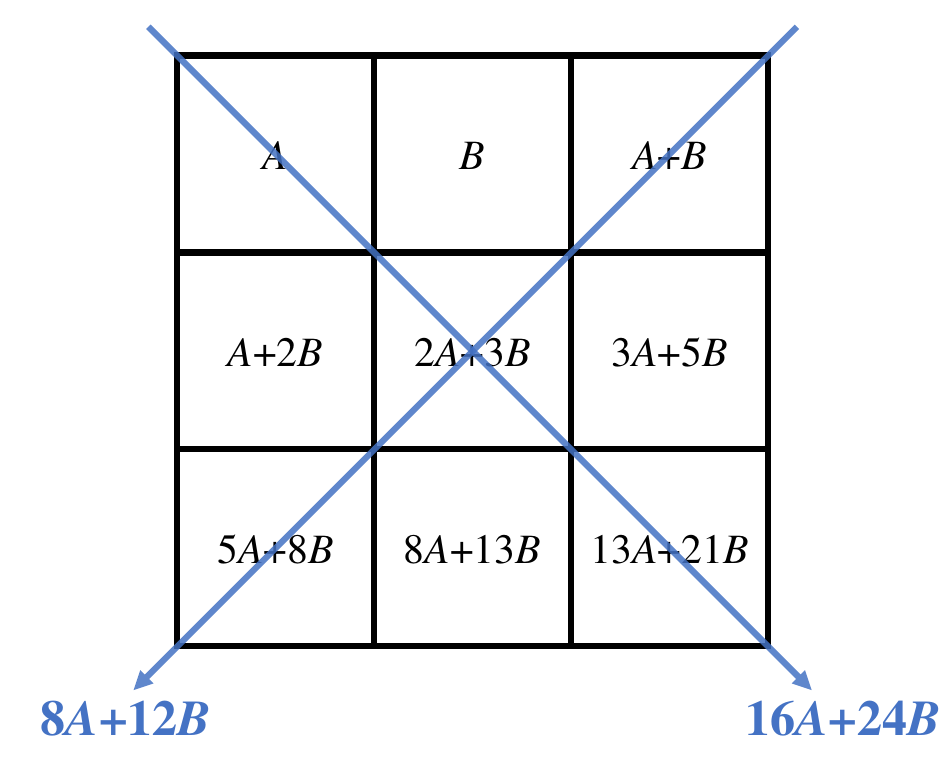}
    \caption{Generalized  Fibonacci  sequence  in  $3\!\times\!3$ grids}
    \label{Generalized Fibonacci sequence in grids}
  \end{figure}

For a $3\times 3$ grid, the following theorem   shows that the ratio of the sum of the main diagonal numbers to the sum of the sub-diagonal numbers is not affected by the initial value of Fibonacci sequence.
\begin{theorem}\label{tr1}
  When the generalized Fibonacci sequence is in a $3\times 3$ grid, the ratio of the sum of the main diagonal numbers to the sum of the sub-diagonal numbers is equal to 2, a constant.
\end{theorem}
{\bf Proof.}
Let the initial values of the generalized Fibonacci sequence be $A, B$. Fill the generalized Fibonacci sequence into the grid, as shown in Figure \ref{Generalized Fibonacci sequence in grids}. Because of the fact that
\[
  \frac{16A+24B}{8A+12B}=2
  \]
   is a constant and is independent of the initial value of the generalized Fibonacci sequence, the proof is completed. $\hfill\blacksquare$

\FloatBarrier
  \FloatBarrier

\end{document}